\input{style/arxiv-ba.cfg}
\documentclass[ba,linksfromyear,preprint]{imsart}
\makeatletter
   \@ifpackageloaded{natbib}{}{\usepackage{natbib}}
\makeatother
\usepackage{xspace}
\newcommand{\itref}[1]{\mbox{\ref{it:#1}}}
\newcommand{\iid}{{independent and identically distributed}\xspace}
\newcommand{\half}{\frac{1}{2}}
\newcommand{\enquote}[1]{``#1''}

\pubyear{2015}
\volume{10}
\issue{2}
\firstpage{517}
\lastpage{521}
\doi{10.1214/15-BA942REJ}

\begin{document}

\begin{frontmatter}
\title{Rejoinder\thanksref{T1}}
%
%
\relateddois{T1}{Main article DOI: \relateddoi[ms=BA942]{Related item:}{10.1214/15-BA942}.}
\runtitle{Rejoinder}

\begin{aug}
\author[a]{\fnms{A.~Philip} \snm{Dawid}\corref{}\ead[label=e1]{apd@statslab.cam.ac.uk}}
\and
\author[b]{\fnms{Monica} \snm{Musio}\ead[label=e2]{mmusio@unica.it}}

\runauthor{A.~P.~Dawid and M.~Musio}

\address[a]{University of Cambridge, \printead{e1}}
\address[b]{University of Cagliari, \printead{e2}}

\end{aug}

\begin{abstract}
We are deeply appreciative of the initiative of the editor, Marina
  Vanucci, in commissioning a discussion of our paper, and extremely
  grateful to all the discussants for their insightful and
  thought-provoking comments.  We respond to the discussions in
  alphabetical order.
\end{abstract}



\end{frontmatter}


\subsection*{Grazian, Masiani and Robert}
Clara Grazian, Ilaria Masiani and Christian Robert (henceforth GMR)
point to a number of potential difficulties in our approach.
\begin{description}
\item[Calibration] We are not sure what GMR mean by the expression
  ``very loosely relates to a loss function.''  A proper scoring rule
  $S(x,Q)$ is very strictly a loss function, where the state is the
  value $x$ of $X$, and the decision is the quoted distribution $Q$
  for $X$.  Moreover (see, for example, \cite{Dawid:1986}), given an
  essentially arbitrary decision problem, with state-space ${\cal X}$,
  decision space ${\cal A}$, and loss function $L(x,a)$, we can define
  $S(x,Q) := L(x, a_Q)$, where $a_Q$ denotes a Bayes act with respect
  to the distribution $Q$ for $X$; and this is readily seen to be a
  proper scoring rule.  That is, essentially every decision problem is
  equivalent to one based on a proper scoring rule.  If you take some
  specified decision problem seriously, you should use the associated
  proper scoring rule.  There is then no problem of calibration.

\item[Dependence on parametrisation] GMR are correct in noting that, if we apply
  a scoring rule after first transforming the state space, we will
  generally get a non-equivalent result (the log-score is essentially
  the only exception to this.)  However, there will be a new scoring
  rule for the transformed problem that is equivalent to the original
  rule for the original problem; see \citet[Section~11]{Parry:2012} for
  how a homogeneous score such as that of Hyv\"arinen transforms.  We
  cannot give any definitive guidance on how to choose an appropriate
  transformation, though Example~11.1 of the above-mentioned paper
  suggests that some consideration of boundary conditions may be
  relevant.

\item[Dependence on dominating measure] This is not the case: when
  constructing the Hyv\"arinen (or other homogeneous) score, the
  formula is to be applied to the density with respect to Lebesgue
  measure.
\item[Arbitrariness] There is indeed a very wide variety of
  homogeneous proper scoring rules, any one of which will achieve our
  aim of eliminating the problematic normalising constant.  At this
  point we can do little more than reiterate what we said towards the
  end of Section~3 of our paper.
\item[Consistency] Whether or not a person, lay or otherwise, finds
  consistency a compelling desideratum is probably a very personal
  matter.  We do find it so.  In a related point, we do not see why,
  in their first paragraph, GMR dismiss the implications of our
  expansion (4) so uncritically.  Indeed, the near identity of the red
  lines in the four subplots of their own Figure~2, which correspond
  to very different prior variances, lends support to our conclusion,
  from (4), that ``the dependence of the Bayes factor on the
  within-model prior specifications is typically negligible.''
\end{description}

GMR correctly point out that there are continuous distributions, such
as the Laplace distribution, to which we cannot apply the Hyv\"arinen
(or other homogeneous) score.  This point deserves further attention.
But for discrete models there is a different class of homogeneous
proper scoring rules that are appropriate and can be used to the same
end of eliminating the normalising constant; see
\cite{apd/sll/mp:2012}.

GMR's simulation studies are interesting.  In contrast to our own
analysis, they appear to show consistency of model selection based on
the multivariate version of the Hyv\"arinen score.  We should not
complain if our method behaves even better than expected, but we
confess we find this puzzling.  We must also take issue with their
assertion that ``the log proper scoring rule tends to infinity
[approximately four times] more slowly than the Bayes factor or than
the likelihood ratio.''  It is simply not appropriate to compare
absolute values across different scoring rules, since each can be
rescaled by an arbitrary positive factor without any consequence for
model comparison.

GMR point to the alternative approach of \cite{kamary}.  However, it
seems to us that the part of that paper that relates to handling
improper priors could just as readily be applied directly to the Bayes
factor.  For example, if we are comparing two location models, we might
use the identical improper prior (with the identical value for its
arbitrary scale factor) for the location parameter in both.  Then this
scale factor will cancel out in the Bayes factor, so leading to an
unambiguous answer.  But in any case, this approach is not available
unless there are parameters in common between all the models being
compared.  Our own approach has no such constraint.

\subsection*{Hans and Perrugia}
Christopher Hans and Mario Perrugia (HP) only consider ``models''
without any unknown parameters, so do not directly address our main
concern, which was to devise methods for comparing parametric models
having possibly improper prior distributions.

They focus on two main issues:
\begin{enumerate}
\item \label{it:comp} Comparisons between the Hyv\"arinen score and
  the log-score.
\item \label{it:robust} Robustness to outliers.
\end{enumerate}

With regard to point~\itref{comp}, HP consider in particular cases
where the two scores are linearly related.  While we fail to see why
this property should be of any fundamental importance (and will pass
up their invitation to characterise it), it is worthy of some
attention.  We do note, however, that, in their analysis of a general
covariance stationary Gaussian process, HP err when they say
``$\sigma^2_{P_i}$ and $\sigma^2_{Q_i}$ are constant in $i$.''  Recall
that $\sigma^2_{P_i}$ is not the unconditional variance of $X_i$
under $P$, but its conditional variance, given $(X_1,\ldots,X_{i-1})$.
Their asserted constancy property will hold for an AR($p$) process
only for $i > p$; while for a general process it will fail, although a
limiting value will typically exist.

HP's specific applications do have this constancy property (at least
for $i>1$).  In the case they consider of different means and equal
variances, the Hyv\"arinen incremental delta score is just a constant
multiple of that for the log-score, and this property extends to the
cumulative scores.  Since an overall positive scale factor is
irrelevant, the two scores are essentially equivalent in this case.

For the other case HP consider, of equal means and different
variances, even after rescaling the incremental delta scores will
differ by an additive constant, $c$ say.  The cumulative scores, to
time $n$, will thus differ by $nc$, which tends to infinity---an
effect that might seem to jeopardise the consistency analysis in our
paper.  However, the following analysis shows that this is not so.
Using HP's formulae, and setting $\xi = \tau_P^2/\tau_Q^2$, consider
first the log-score.  The incremental delta log-scores are, under $P$,
\iid, with expectation $\half(\xi- 1 -\log \xi) > 0$ and finite
variance, so that the difference between the cumulative prequential
score for $Q$ and that for $P$ tends to infinity almost surely---so
favouring the true model $P$.  Likewise $Q$ will be favoured when it
is true.  Now consider the Hyv\"arinen score.  Again the incremental
delta log-scores under $P$ are \iid with finite variance, now with
expectation $\tau_q^{-2}(\xi + \xi^{-1} -2) > 0$; so once again, the
true model is consistently favoured.

HP ask whether there is any principled reason for applying the cut-off
value $0$ to the difference in prequential scores.  Well, it seems
natural to us to choose the model whose predictions have performed
best so far, so indicating that this might continue into the
future---although, as the advertisers of financial products are
obliged to point out, past success cannot be taken as an infallible
guide to future performance.  We further note the essential
equivalence of this recipe to the machine learning technique of
``empirical risk minimisation'' in Statistical Learning Theory, which
has developed an extensive theory, extending well beyond the case of
parametric models, characterising when this will be effective; see
\cite{rakhlin15sequential,sasha:chapter} for application to the
general case of dependent sequential observations.

In any case, should one wish to use a cut-off different from $0$,
there is no impediment to doing so---this would not affect the
consistency properties we have investigated, which only rely on the
difference of cumulative scores tending to infinity.  How the choice
of cut-off could relate to differential prior probabilities and
utilities is a topic that deserves further consideration.

Turning to HP's point~\itref{robust}, their simulations appear to show
that the Hyv\"arinen score is less robust to additive outliers than
the log-score (though we note that in their example the outlier only
affects 2 of the 100 summands of the overall score.)  Issues of the
robustness of minimum score inference have been considered by
\cite{apd/mm/lv}, where it is shown that (in an estimation context)
certain proper scoring rules do enjoy good robustness properties
(generally better than straightforward likelihood).  However, these do
not include the Hyv\"arinen score or other homogeneous scores.  Thus
there may indeed be a conflict between the aim of our current paper,
which is to overcome problems associated with improper distributions,
and the very different aim of protecting against outliers.

\subsection*{Katzfuss and Bhattacharya}
Matthias Katzfuss and Anirban Bhattacharya (KB) are particularly
concerned with the question of whether our approach can be tweaked to
yield a ``pseudo-Bayes factor'', where a general score takes the place
of log-likelihood.  While it would be very nice if this were so, we
are a little dubious.  As KB point out, there are serious problems
related to the arbitrary scaling of a general score.  These are
compounded when, as for the homogenous cases we consider, the score is
a dimensioned quantity.  Thus if the basic observable $X$ has the
dimension of length, $L$, then the Hyv\"arinen score has dimension
$L^{-2}$, so any scale factor, such as $\lambda$ in their (1.1) or
(3.1), would have to have dimension $L^2$.  Otherwise put, whether we
are measuring $X$ in nanometers or in parsecs will affect the absolute
value of the score (though not the comparisons that form the basis of
our method).

There is no reason why our method should not be used to compare a
finite number of models, rather than just 2.  However, when the number
is countably infinite, or grows with sample size, even
likelihood-based model selection can fail to be consistent.  In that
case the problem can sometimes be solved by regularisation,
essentially equivalent to introducing prior probabilities over models
and selecting on the basis of the posterior model probabilities.
Perhaps some analogue of this device might work for more general
proper scoring rules.





\end{document}